\documentclass[11pt]{article}
\usepackage{amssymb,latexsym}
\usepackage{epsfig}
\usepackage{eufrak}
\usepackage{amsmath}
\usepackage{mathrsfs}
\usepackage{color}

\newtheorem{theorem}{Theorem}
\newtheorem{lemma}{Lemma}

\newcommand{\be}{\begin{equation}}
\newcommand{\ee}{\end{equation}}
\newcommand{\bee}{\begin{eqnarray*}}
\newcommand{\eee}{\end{eqnarray*}}
\newcommand{\bel}{\begin{eqnarray}}
\newcommand{\eel}{\end{eqnarray}}
\newcommand{\bec}{\begin{cases}}
\newcommand{\eec}{\end{cases}}
\newcommand{\bem}{\begin{bmatrix}}
\newcommand{\eem}{\end{bmatrix}}

\newcommand{\la}{\label}
\newcommand{\li}{\left}
\newcommand{\ri}{\right}

\newcommand{\lc}{\lceil}
\newcommand{\rc}{\rceil}
\newcommand{\lf}{\lfloor}
\newcommand{\rf}{\rfloor}

\newcommand{\vep}{\varepsilon}

\newcommand{\de}{\delta}
\newcommand{\De}{\Delta}

\newcommand{\ro}{\rho}

\newcommand{\f}{\frac}

\newcommand{\qu}{\quad}
\newcommand{\qqu}{\qquad}

\newcommand{\mscr}{\mathscr}

\newcommand{\bb}{\mathbb}

\newcommand{\mrm}{\mathrm}

\newcommand{\LRA}{\Longleftrightarrow}
\newcommand{\sh}{\slash}

\newcommand{\tx}{\text}

\newcommand{\bed}{\begin{description}}
\newcommand{\eed}{\end{description}}
\newcommand{\bei}{\begin{itemize}}
\newcommand{\eei}{\end{itemize}}
\newcommand{\ben}{\begin{enumerate}}
\newcommand{\een}{\end{enumerate}}

\newcommand{\beL}{\begin{lemma}}
\newcommand{\eeL}{\end{lemma}}
\newcommand{\beT}{\begin{theorem}}
\newcommand{\eeT}{\end{theorem}}
\newcommand{\sect}{\section}

\newcommand{\bpf}{\begin{pf}}
\newcommand{\epf}{\end{pf}}
\newcommand{\bsk}{\bigskip}
\newcommand{\bi}{\binom}

\newcommand{\pfbox}{\hfill\mbox{$\Box$}}

\newenvironment{pf}{\paragraph*{Proof{\rm.}}}{\pfbox\bigskip}

\begin{document}

\title{{\bf Exact Computation of Minimum Sample Size for Estimating Proportion of Finite Population}
\thanks{The author had been previously working with Louisiana State University at Baton Rouge, LA 70803, USA,  and is now with
Department of Electrical Engineering, Southern University and A\&M College, Baton Rouge, LA 70813, USA; Email:
chenxinjia@gmail.com}}

\author{Xinjia Chen}

\date{July 2007}

\maketitle

\begin{abstract}

In this paper, we develop an exact method for the determination of
the minimum sample size for estimating the proportion of a finite
population with prescribed margin of error and confidence level.  By
characterizing the behavior of the coverage probability with respect
to the proportion, we show that the computational complexity can be
significantly reduced and bounded regardless of population size.

\end{abstract}

\sect{Introduction}

The estimation of the proportion of a finite population is a basic and very important problem in probability and
statistics \cite{Desu, Thompson}. The problem is formulated as follows.

Consider a finite population of $N$ units, among which there are $M$ units having a certain attribute. It is a
frequent problem to estimate the proportion $\f{M}{N}$ by sampling without replacement.  Let $n$ be the sample
size and $\mathbf{k}$ be the number of units that found to carry the attribute. The estimate of the proportion
is taken as $\f{\mathbf{k}}{n}$.  A crucial question in the estimation is as follows:

\bsk

{\it Given the knowledge that $M$ belongs to interval $[L, U]$, what
is the minimum sample size $n$ that guarantees the difference
between $\f{\mathbf{k}}{n}$ and $\f{M}{N}$ be bounded within some
prescribed margin of error with a confidence level higher than a
prescribed value?}

\bsk

Conventionally, the exact method requires evaluation of  the
coverage probability for all values of $M$ in $[L, U]$ for sample
sizes incrementing from $2$ to a number large enough.  Since the
range of interval $[L, U]$ can be as wide as $[0, N]$,  the number
of evaluations of coverage probability can be very large if the
population size $N$ is large.  The main contribution of this paper
is to provide exact method for the computation of minimum sample
size such that the total number of evaluations of coverage
probability can be significantly reduced and bounded regardless of
the population size $N$. Specially, we demonstrate that a small
subset of the integers in interval $[L, U]$ needs to be evaluated.

The paper is organized as follows. In Section 2, the techniques for computing the minimum sample size is
developed with the margin of error taken as a bound of absolute error. In Section 3, we derive corresponding
sample size method by using relative error bound as the margin of error.  In Section 4, we develop techniques
for computing minimum sample size with a mixed error criterion.  Section 5 is the conclusion. The proofs are
given in Appendices.

Throughout this paper, we shall use the following notations. The set of integers is denoted by $\bb{Z}$.  The
ceiling function and floor function are denoted respectively by $\lc . \rc$ and $\lf . \rf$ (i.e., $\lc x \rc$
represents the smallest integer no less than $x$; $\lf x \rf$ represents the largest integer no greater than
$x$).  For non-negative integer $m$, the combinatoric function $\bi{m}{z}$ with respect to integer $z$ means
\[
\bi{m}{z} = \bec \f{ m! } { z! (m- z)! } & \tx{for} \; 0 \leq z \leq m,\\
0 & \tx{for} \; z < 0 \; \tx{or} \; z > m. \eec
\]
We denote
\[
S(n,k,l,M,N) = \sum_{i=k}^l \f{ \bi{M}{i} \bi{N-M}{n- i} } { \bi{N}{n} }.
\]
The notation ``$\LRA$'' means ``if and only if''. The other notations will be made clear as we proceed.

 \sect{Control of Absolute Error}

 For $\vep, \; \de \in (0,1)$, it is desirable in many situations to find the
minimum sample size $n$ such that
\[
\Pr \li \{ \li | \f{\mathbf{k}}{n} - \f{M}{N}  \ri | < \vep \ri \}
> 1 -  \de
\]
for any $M$ in interval $[L, U]$.  Here the interval $[L, U]$ is introduced to take into account the knowledge
of $M$.  If no information about $M$ is available, $[L, U]$ is taken as $[0, N]$.  The quantity $\Pr \li \{ \li
| \f{\mathbf{k}}{n} - \f{M}{N}  \ri | < \vep \ri \}$ is referred to as the coverage probability. An essential
step to find the minimum sample size is to determine whether a fixed sample size $n$ is large enough to ensure
that the coverage probability is above $1 - \de$ for any $M$ in $[L, U]$.  By the conventional method, for a
fixed sample size $n$, the total number of evaluations of the coverage probability is $U - L + 1$.  The
computation can be high for large population size $N$.  Interestingly, we discovered that the number of
evaluations of the coverage probability can be significantly reduced by taking advantage of the coverage
property as characterized by Theorem \ref{thm_abs} at below.

 \beT
 \la{thm_abs}  Let $N, \; n$ and $\vep$ be fixed. Let $L$ and $U$
 be two integers such that $0 \leq L \leq U \leq N$.
 Suppose that $L \leq M \leq U$.  Then,  the
minimum of $\Pr \li \{ \li | \f{\mathbf{k}}{n} - \f{M}{N}  \ri | < \vep \ri \}$ with respect to $M$  is attained
at {\small $\{L, U \} \cup \li \{ \li \lf N \li ( \f{k}{n} - \vep \ri ) \ri \rf  \in (L, U): k \in \bb{Z} \ri \}
\bigcup \li \{ \li \lc N \li ( \f{k}{n} + \vep \ri ) \ri \rc \in (L, U): k \in \bb{Z} \ri \}$, } which has less
than $2 n \li( \f{ U - L - 1} { N } \ri ) + 4$ elements.
 \eeT

 See Appendix A for a proof.

\bsk

By the fact of symmetry that $\Pr \li \{ \li | \f{\mathbf{k}}{n} - \f{M}{N} \ri | \ri \} = \Pr \li \{ \li | \f{n
- \mathbf{k}}{n} - \f{N - M}{N} \ri | \ri \}$, we can restrict $M$ to be no larger than $\lc \f{N}{2} \rc$.
Hence, without loss of generality, we can assume that $0 \leq L < U \leq \lc \f{ N}{ 2 } \rc$.  Specially, for
$L = 0, \; U = \lc \f{ N}{ 2 } \rc$,  the total number of evaluations of coverage probability is less than
\[
2 n \li ( \f{\lc \f{ N}{ 2 } \rc -1} { N } \ri ) + 3
 <  2 n \li ( \f{\f{ N}{ 2 } } { N } \ri ) + 3 =  n + 3
\]
since the coverage probability for $L = 0$ is $1$.  This means that,  in the situation that no information about
$M$ is available, the total number of evaluations of coverage probability is at most $\boxed{n + 2}$, which is
independent of the population size $N$. \bsk

\sect{Control of Relative Error}

For $\vep, \; \de \in (0,1)$, it is interesting to find the minimum sample size such that
\[
\Pr \li \{ \li | \f{\mathbf{k}}{n} - \f{M}{N}  \ri | < \vep \f{M}{N}
\ri \}
> 1 -  \de
\]
for any $M$ in the interval $[L, U]$.  For the purpose of reducing the number of evaluations of the coverage
probability, we have
 \beT
 \la{thm_rev}  Let $N, \; n$ and $\vep$ be fixed.
 Let $L$ and $U$ be two integers such that $0 \leq L \leq U \leq N$.
 Suppose that $L \leq M \leq U$. Then,  the
minimum of $\Pr \li \{ \li | \f{\mathbf{k}}{n} - \f{M}{N}  \ri | <
\vep \f{M}{N} \ri \}$ with respect to $M$  is attained at {\small
\[ \{L, U \} \cup \li \{ \li \lf  \f{N k}{(1 + \vep) n} \ri \rf
\in (L, U):  k \in \bb{Z} \ri \} \bigcup \li \{ \li \lc \f{N k}{(1 -
\vep)n} \ri \rc \in (L, U): k \in \bb{Z} \ri \},
\]}
which has less than $2 n \li( \f{ U - L - 1} { N } \ri ) + 4$
elements.
 \eeT

We omit the proof of Theorem \ref{thm_rev} because of its similarity to that of Theorem \ref{thm_abs}.

 \sect{Control of Absolute Error or Relative Error}

Let $\vep_a \in (0,1)$ and $\vep_r \in (0,1)$ be respectively the margins of absolute and relative error.  Let
$\de \in (0,1)$, it is frequently useful to find the minimum sample size such that
\[
\Pr \li \{ \li | \f{\mathbf{k}}{n} - \f{M}{N}  \ri | < \vep_a \;\;
\tx{or} \;\; \li | \f{\mathbf{k}}{n} - \f{M}{N}  \ri | < \vep_r
\f{M}{N} \ri \}
> 1 -  \de
\]
for any $M$ in the interval $[L, U]$.  To reduce the computational complexity, we have

 \beT
 \la{thm_abs_rev}  Let $N, \; n, \; \vep_a$ and $\vep_r$ be fixed.
Let $L$ and $U$ be two integers such that $0 \leq L < \f{N
\vep_a}{\vep_r} < U \leq N$.  Suppose that $L \leq M \leq U$.
 Then, the minimum of $\Pr \li \{ \li | \f{\mathbf{k}}{n} -
\f{M}{N}  \ri | < \vep_a \; \mrm{or} \;  \li | \f{\mathbf{k}}{n} -
\f{M}{N}  \ri | < \vep_r \f{M}{N} \ri \}$ with respect to $M$  is
attained at {\small \bee \li \{L, U, \li \lf \f{N \vep_a}{\vep_r}
\ri \rf, 1 + \li \lf \f{N \vep_a}{\vep_r} \ri \rf \ri \} & \bigcup
& \li \{ \li \lf N \li (  \f{k}{n} - \vep_a \ri ) \ri \rf \in \li
(L, \li \lf \f{N \vep_a}{\vep_r} \ri \rf \ri): k \in
\bb{Z} \ri \}\\
&    \bigcup & \li \{ \li \lc N \li ( \f{k}{n} + \vep_a \ri ) \ri
\rc \in \li (L, \li \lf \f{N
\vep_a}{\vep_r} \ri \rf \ri ): k \in \bb{Z} \ri \}\\
&     \bigcup  & \li \{ \li \lf  \f{N k}{(1 + \vep_r) n} \ri \rf
\in \li (1 + \li \lf \f{N \vep_a}{\vep_r} \ri \rf,  U \ri ):  k
\in \bb{Z} \ri
\}\\
&   \bigcup & \li \{ \li \lc \f{N k}{(1 - \vep_r)n} \ri \rc \in
\li ( 1 + \li \lf \f{N \vep_a}{\vep_r} \ri \rf, U \ri ): k \in
\bb{Z} \ri \}, \eee} which has less than $2 n \li( \f{ U - L - 3}
{ N } \ri ) + 8$ elements.
 \eeT

It should be noted that Theorem \ref{thm_abs_rev} can be shown by applying Theorem \ref{thm_abs} and Theorem
\ref{thm_rev} with the observation that {\small \[ \Pr \li \{ \li | \f{\mathbf{k}}{n} - \f{M}{N}  \ri | < \vep_a
\; \; \mrm{or} \;\; \li | \f{\mathbf{k}}{n} - \f{M}{N}  \ri | < \vep_r \f{M}{N} \ri \} = \bec \Pr \li \{ \li |
\f{\mathbf{k}}{n} - \f{M}{N} \ri | < \vep_a  \ri \} & \tx{for} \; M \in \li [ L,
\li \lf \f{N \vep_a}{\vep_r} \ri \rf \ri ], \\
\Pr \li \{  \li | \f{\mathbf{k}}{n} - \f{M}{N}  \ri | < \vep_r
\f{M}{N}  \ri \} & \tx{for} \; M \in \li [ 1 + \li \lf \f{N
\vep_a}{\vep_r} \ri \rf, U \ri ]. \eec
\]}

Finally, we would like to point out that similar characteristics of coverage probability can be shown for other
populations such as Bernoulli population and Poisson population, which allows for the exact computation of
minimum sample size. For details, see our recent papers \cite{Chen, Chen2}.

\section{Conclusion}

In this paper, we develop an exact method for the computation of the
minimum sample size for estimation of the proportion of finite
population.  The method is much more efficient than previously
possible.  The efficiency improvement is due to the interesting
discovery of the characteristics of the coverage probability. Such
characteristics reveals a new aspect of the hyper-geometrical
distribution.

\appendix

\sect{Proof of Theorem \ref{thm_abs}}

Define
\[
C(M) = \Pr \li \{ \li | \f{\mathbf{k}}{n} - \f{M}{N}  \ri | < \vep
\ri \} = \Pr \li \{ g(M) \leq \mathbf{k} \leq  h(M)  \ri \}
\]
where
\[
g(M) = \li \lf n \li (  \f{M}{N} - \vep \ri ) \ri \rf  + 1, \qqu
h(M) = \li \lc n \li (  \f{M}{N} + \vep \ri ) \ri \rc - 1.
\]
It should be noted that $C(M), \; g(M)$ and $h(M)$ are actually
multivariate functions of $M, \; N, \; \vep$ and $n$.  For
simplicity of notations, we drop the arguments $n, \; N$ and
$\vep$ throughout the proof of Theorem \ref{thm_abs}.

\beL \label{tau increase} Let $ 0 \leq M < M+1 \leq N$. Define $T(k,
M, N, n) = \li. {M \choose k} {N-M-1 \choose n-k-1} \ri \slash {N
\choose n}$.  Then, $S (n,0,k,M, N) - S (n,0,k,M+1, N) = T(k, M, N,
n)$ for any integer $k$. \eeL

\bpf

We first show the equation for $0 \leq k \leq M$.  We perform
induction on $k$. For $k = 0$, we
have \bel S (n,0,k,M, N) - S (n,0,k,M+1, N) & = & S (n,0,0,M, N) - S (n,0,0,M+1, N) \nonumber\\
& = &   \frac{ {M \choose 0} {N-M \choose n} }
 { {N \choose n} } -  \frac{ {M + 1 \choose 0} {N-M-1 \choose n} }
 { {N \choose n} } \nonumber \\
 & =  & \frac{ {N-M-1 \choose n - 1} }
 { {N \choose n} } \la{comb}\\
 & =  & \frac{ {M \choose 0} {N-M-1 \choose n-0-1} }
 { {N \choose n} } = T(0, M, N, n) \nonumber, \eel where
 (\ref{comb}) follows from the fact that, for non-negative integer
 $m$, \be \la{basic} \bi{m + 1}{z+1} = \bi{m}{z} +
\bi{m}{z+1} \ee for any integer $z$.

Now suppose the lemma is true for $k-1$ with $1 \leq k \leq M$,
i.e.,
 \[
S (n,0,k-1,M, N) - S (n,0,k-1,M+1, N)  = \frac{ {M \choose k-1}
{N-M-1 \choose n-k} }
 { {N \choose n} }.
 \]
Then, {\small \bel  S (n,0,k,M, N) - S (n,0,k,M+1, N)
 & = &  S (n,0,k-1,M, N) - S (n,0,k-1,M+1, N) \nonumber\\
 &   &  + \frac{ {M \choose k} {N-M
\choose n - k} } { {N \choose n} } -  \frac{ {M + 1 \choose k}
{N-M-1 \choose n - k} }
 { {N \choose n} } \nonumber\\
 & = & \frac{ {M \choose k-1} {N-M-1 \choose n-k} }
 { {N \choose n} } +  \frac{ {M \choose k} {N-M
\choose n - k} } { {N \choose n} } -  \frac{ {M + 1 \choose k}
{N-M-1 \choose n - k} }
 { {N \choose n} } \nonumber\\
 &  =  & \frac{ {M \choose k} {N-M
\choose n - k} } { {N \choose n} } - \li [ \frac{ {M + 1 \choose
k} {N-M-1 \choose n - k} }
 { {N \choose n} } - \frac{ {M \choose k-1} {N-M-1 \choose n-k} }
 { {N \choose n} } \ri ] \nonumber\\
& = & \frac{ {M \choose k} {N-M \choose n - k} } { {N \choose n} }
- \frac{ {M \choose k} {N-M - 1 \choose n -
k} } { {N \choose n} } \la{eqa}\\
& = & \frac{ {M \choose k} {N-M - 1 \choose n - k - 1} } { {N
\choose n} } \la{eqb} \eel} where (\ref{eqa}) and (\ref{eqb})
follows from (\ref{basic}).  Therefore, we have shown the lemma
for $0 \leq k \leq M$.

For $k > M$, we have $S (n,0,k,M, N) = S (n,0,k,M+1, N) = 1$ and
$T(k, M, N, n) = 0$.  For $k < 0$, we have $S (n,0,k,M, N) = S
(n,0,k,M+1, N) = 0$ and $T(k, M, N, n) = 0$.  Thus, the lemma is
true for any integer $k$.

\epf

\beL \la{dou} Let $1 \leq M \leq N$ and $k \leq l$. Then,
\[ S(n, k, l, M, N) - S(n, k, l, M - 1, N) = T(k - 1,
M-1, N, n) - T(l, M - 1, N, n).
\]
\eeL

\bpf

To show the lemma, it suffices to consider $6$ cases as follows.

Case (i): $0 < n < k \leq l$. In this case, $S(n, k, l, M, N) = S(n,
k, l, M - 1, N) = 0$ and $T(k - 1, M-1, N, n) = T(l, M - 1, N, n) =
0$.

Case (ii): $k \leq l < 0 < n$. In this case, $S(n, k, l, M, N) =
S(n, k, l, M - 1, N) = 0$ and $T(k - 1, M-1, N, n) = T(l, M - 1, N,
n) = 0$.

Case (iii): $k \leq 0 < n \leq l$.  In this case, $S(n, k, l, M, N)
= S(n, k, l, M - 1, N) = 1$ and $T(k - 1, M-1, N, n) = T(l, M - 1,
N, n) = 0$.

Case (iv): $k \leq 0 \leq l < n$.  In this case, $T(k - 1, M-1, N,
n) = 0$ and, by Lemma \ref{tau increase},  {\small \bee S(n, k, l,
M, N) - S(n, k, l, M - 1, N)
& = & [S(n, 0, l, M, N) - S(n, 0, l, M - 1, N)]\\
& = & T(k - 1, M-1, N, n) - T(l, M - 1, N, n).
 \eee}

Case (v): $0 < k \leq n \leq l$.  In this case, $T(l, M-1, N, n) =
0$ and, by Lemma \ref{tau increase}, {\small  \bee S(n, k, l, M, N)
- S(n, k, l, M - 1, N)
& = & [S(n, 0, k-1, M - 1, N) - S(n, 0, k-1, M, N)]\\
& = & T(k - 1, M-1, N, n) - T(l, M - 1, N, n).
 \eee}

Case (vi): $0 < k \leq l < n$.  In this case, by Lemma \ref{tau
increase},  {\small \bee S(n, k, l, M, N) - S(n, k, l, M - 1, N)
& = & [S(n, 0, l, M, N) - S(n, 0, k - 1, M, N)]\\
&   &  - [S(n, 0, l, M - 1, N) - S(n, 0, k - 1, M - 1, N)]\\
& = & [S(n, 0, l, M, N) - S(n, 0, l, M - 1, N)]\\
&   &  - [S(n, 0, k - 1, M, N) - S(n, 0, k - 1, M - 1, N)]\\
& = & T(k - 1, M-1, N, n) - T(l, M - 1, N, n).
 \eee}

\epf

\beL \la{lem22} Let $l \geq 0$ and $k < n$.  Then, $\li \lf \f{ n M}
{ N+1} \ri \rf \geq l$ for $M \geq 1 + \li \lf \f{ N l } { n - 1 }
\ri \rf$, and  $\li \lf \f{  n M} { N+1} \ri \rf \leq k - 1$ for $M
\leq 1 + \li \lf \f{ N (k - 1) } { n - 1 } \ri \rf$.

\eeL

\bpf

 To show the first part of the lemma, observe that $(N+1 - n) l \geq 0$,
 by which we can show $\f{n N l}{n-1}
\geq (N+1) l$. Hence, $n \li ( 1 + \li \lf \f{ N l } { n - 1 } \ri
\rf \ri )
> \f{n N l}{n-1} \geq (N+1)
l$.  That is, $ \f{n} { N+1} \li ( 1 + \li \lf \f{ N l } { n - 1 }
\ri \rf \ri ) >  l$.  It follows that $\li \lf \f{n}{ N+1} \li ( 1 +
\li \lf \f{ N l } { n - 1 } \ri \rf \ri ) \ri \rf \geq l$. Since the
floor function is non-decreasing, we have $\li \lf \f{ n M} { N+1}
\ri \rf \geq l$ for $M \geq 1 + \li \lf \f{ N l } { n - 1 } \ri
\rf$.

To prove the second part of the lemma, note that $(N + 1 - n) (n -
k) > 0$, from which we can deduce $1 + \f{ N(k - 1) } { n - 1 } <
\f{ (N +1) k } { n }$.  Hence, $1 + \li \lf \f{ N(k - 1) } { n - 1 }
\ri \rf < \f{ (N +1) k  } { n }$, i.e., $\f{ n } {N + 1} \li ( 1 +
\li \lf \f{ N(k - 1) } { n - 1 } \ri \rf \ri ) < k$, leading to $\li
\lf \f{ n } {N + 1} \li ( 1 + \li \lf \f{ N(k - 1) } { n - 1 } \ri
\rf \ri ) \ri \rf \leq k - 1$. Since the floor function is
non-decreasing, we have $\li \lf \f{ n M} { N+1} \ri \rf \leq k - 1$
for $M \leq 1 + \li \lf \f{ N (k - 1) } { n - 1 } \ri \rf$.

\epf

\beL \la{monotone}
 Let $0 \leq r \leq n$.  Then, the following statements hold true.

(I) \[ T(r - 1, M - 1, N, n) \leq T(r, M - 1, N, n) \qu \mrm{for}
\qu 1 \leq r \leq \li \lf \f{ n M} { N+1} \ri \rf;
\]
\[ T(r + 1, M - 1, N, n) \leq T(r, M - 1, N, n) \qu \mrm{for}
\qu  \li \lf \f{ n M} { N+1} \ri \rf \leq r \leq n - 1.
\]

(II)
\[
T(r, M - 2, N, n) \leq T(r, M - 1, N, n) \qu \mrm{for} \qu 1 < M
\leq 1 + \li \lf \f{ N r } { n - 1 } \ri \rf;
\]
\[
T(r, M, N, n) \leq T(r, M - 1, N, n) \qu \mrm{for} \qu  1 + \li
\lf \f{ N r } { n - 1 } \ri \rf \leq M < N.
\]

\eeL

\bpf

To show statement (I), note that $T(r, M -1, N, n) = 0$ for
$\min(M - 1, n-1) < r \leq n$.  Our calculation shows that
 \[
 \f{ T(r - 1, M -1, N, n) } { T(r, M -1, N, n)
} = \f{r}{ M - r} \f{ N - M + 1 - (n -r) } { n - r } \leq 1 \qu
\mrm{for} \qu 1 \leq  r \leq \f{ n M } { N + 1 } \]
 and
 \[
 \f{ T(r - 1, M-1, N, n) }
{ T(r, M -1, N, n) } > 1 \qu \mrm{for} \qu \f{ n M } { N + 1 } < r
\leq \min(M - 1, n-1).
\]

To show statement (II), note that $T(r, M - 1, N, n) = 0$ for $1
\leq M < r + 1$, and $T(r, M - 1, N, n) \geq T(r, M - 2, N, n) = 0$
for $M = r + 1$.  Direct computation shows that {\small \[ \f{ T(r,
M -1, N, n) } { T(r, M - 2, N, n) } = \f{ M -1 } {M - 1 - r}
 \f{ N - M  + 2 - (n -r) } {N - M + 1} \geq 1 \qu \mrm{ for} \qu r + 1 < M  \leq 1 +  \f{ N r } { n - 1
 },
 \]}
  and
  \[
  \f{ T(r, M - 1, N, n) } { T(r, M - 2, N,
n) } < 1 \qu \mrm{ for } \qu  1 + \f{ N r } { n - 1 } < M \leq N.
\]

\epf

\beL \la{lem77}
Let $k \leq l$ and $0 \leq L \leq U \leq N$.  Then,
\[
\min_{M \in [L, \; U]}  S(n, k, l, M, N) = \min \{ S(n, k, l, L,
N), \; S(n, k, l, U, N) \}.
\]
\eeL

\bpf

To show the lemma, it suffices to consider $6$ cases as follows.

Case (i): $0 < n < k \leq l$. In this case, $S(n, k, l, M, N) = 0$
for any $M \in [L, U]$.

Case (ii): $k \leq l < 0 < n$. In this case, $S(n, k, l, M, N) = 0$
for any $M \in [L, U]$.

Case (iii): $k \leq 0 < n \leq l$. In this case, $S(n, k, l, M, N) =
1$ for any $M \in [L, U]$.

Case (iv): $k \leq 0 \leq l < n$.  In this case, $S(n, k, l, M, N) =
S(n, 0, l, M, N)$ is non-increasing with respect to $M \in [L, U]$
as can be seen from Lemma \ref{tau increase}.

Case (v): $0 < k \leq n \leq l$. In this case, $S(n, k, l, M, N) =
1- S(n, 0, k - 1, M, N)$ is non-decreasing with respect to $M \in
[L, U]$ as can be seen from Lemma \ref{tau increase}.

Clearly, the lemma is true for the above five cases.

Case (vi): $0 < k \leq l < n$.  Define {\small $\De(k, l, M, N, n)
= S(n, k, l, M, N) - S(n, k, l, M - 1, N)$}.  By Lemma \ref{dou},
$\De(k, l, M, N, n) = T(k - 1, M-1, N, n) - T(l, M - 1, N, n)$.

Invoking Lemma \ref{lem22}, for $M \geq 1 + \li \lf \f{ N l } { n -
1 } \ri \rf$, we have that $\li \lf \f{  n M} { N+1} \ri \rf \geq l$
and thus, by statement (I) of Lemma \ref{monotone}, $T(r, M - 1, N,
n)$ is non-decreasing with respect to $r \leq l$. Consequently, $T(k
- 1, M -1, N, n) \leq T(l, M -1, N, n)$, leading to $\De(k, l, M, N,
n) \leq 0$ for $M \geq 1 + \li \lf \f{ N l } { n - 1 } \ri \rf$.

Similarly, applying Lemma \ref{lem22}, for $M \leq 1 + \li \lf \f{ N
(k - 1) } { n - 1 } \ri \rf$, we have that $\li \lf \f{  n M} { N+1}
\ri \rf \leq k - 1$ and thus, by statement (I) of Lemma
\ref{monotone}, $T(r, M - 1, N, n)$ is non-increasing with respect
to $r \geq k - 1$. Consequently, $T(k -1, M -1, N, n) \geq T(l, M
-1, N, n)$, leading to $\De(k, l, M, N, n) \geq 0$ for $M \leq 1 +
\li \lf \f{ N (k - 1) } { n - 1 } \ri \rf$.

By statement (II) of Lemma \ref{monotone}, for $1 + \li \lf \f{ N (k
- 1) } { n - 1 } \ri \rf \leq M \leq 1 + \li \lf \f{ N l } { n - 1 }
\ri \rf$, we have that $T(l, M - 1, N, n)$ is non-decreasing with
respect to $M$ and that $T(k - 1, M - 1, N, n)$ is non-increasing
with respect to $M$. It follows that $\De(k, l, M, N, n)$ is
non-increasing with respect to $M$ in this range.  Therefore, there
exists an integer $M^*$ such that $1 + \li \lf \f{ N (k - 1) } { n -
1 } \ri \rf \leq M^* \leq 1 + \li \lf \f{ N l } { n - 1 } \ri \rf$
and that $\De(k, l, M, N, n) \geq 0$ for $0 \leq M \leq M^*$, and
$\De(k, l, M, N, n) \leq 0$ for $M^* \leq M \leq N$.  This implies
that $S(n, k, l, M, N)$ is non-decreasing for $0 \leq M \leq M^*$
and non-increasing for $M^* \leq M \leq N$. The concludes the proof
of the lemma.

\epf

\beL \la{lemh} Let $m$ be an integer.  Then,
\[
 \li \lc n \li (  \f{m}{N} + \vep \ri ) \ri \rc  = \bec k &
\tx{for} \; m = \li \lf N \li (  \f{k}{n} - \vep \ri ) \ri \rf,\\
k + 1 & \tx{for} \; \li \lf N \li (  \f{k}{n} - \vep \ri ) \ri \rf
< m \leq \li \lf N \li (  \f{k + 1}{n} - \vep \ri ) \ri \rf. \eec
\]
 \eeL

\bpf For notational simplicity, let {\small $r = \li \lf N \li (
\f{k}{n} - \vep \ri ) \ri \rf$} and {\small $r^\prime = \li \lf N
\li ( \f{k+1}{n} - \vep \ri ) \ri \rf$}. By the definition of the
floor function, we have {\small $N \li (  \f{k}{n} - \vep \ri ) - 1
< r \leq N \li (  \f{k}{n} - \vep \ri )$}, which can be written as
{\small $k - \f{n}{N} < n \li (  \f{r}{N} + \vep \ri ) \leq k$}. As
a result, {\small $\li \lc n \li (  \f{r}{N} + \vep \ri ) \ri \rc  =
k$} and the lemma is true for $m = r$.  Similarly, {\small $\li \lc
n \li (  \f{r^\prime}{N} + \vep \ri ) \ri \rc = k + 1$} and the
lemma is true for $m = r^\prime$.

Since $r$ is an integer, we have {\small $N \li (  \f{k}{n} - \vep
\ri ) < r+1 \leq m \leq r^\prime -1 \leq N \li ( \f{k+1}{n} - \vep
\ri ) - 1$} for $r < m < r^\prime$.  Hence, {\small $\f{k}{n} <
\f{m}{N} + \vep \leq \f{k+1}{n} - \f{1}{N}$}.  That is, {\small $k <
n \li ( \f{m}{N} + \vep \ri ) \leq k + 1 - \f{n}{N}$}, which implies
that {\small $\li \lc n \li ( \f{m}{N} + \vep \ri ) \ri \rc = k +
1$} for $r < m < r^\prime$. The proof of the lemma is thus
completed.

\epf

\beL \la{lemg} Let $m$ be an integer.  Then,
\[
\li \lf n \li (  \f{m}{N} - \vep \ri ) \ri \rf  = \bec k &
\tx{for} \; \li \lc N \li (  \f{k}{n} + \vep \ri ) \ri \rc \leq m
< \li \lc N \li (  \f{k + 1}{n} +\vep \ri ) \ri \rc,\\
k + 1 & \tx{for} \; m = \li \lc N \li (  \f{k+1}{n} + \vep \ri )
\ri \rc. \eec
\]
 \eeL

\bpf For notational simplicity, let {\small $r = \li \lc N \li (
\f{k}{n} + \vep \ri ) \ri \rc$} and {\small $r^\prime = \li \lc N
\li ( \f{k+1}{n} + \vep \ri ) \ri \rc$}.  By the definition of the
ceiling function, we have {\small $N \li (  \f{k}{n} + \vep \ri )
\leq r < N \li (  \f{k}{n} + \vep \ri ) + 1$}, which can be written
as {\small $k \leq  n \li (  \f{r}{N} - \vep \ri ) < k + \f{n}{N}$}.
Hence, {\small $\li \lf n \li (  \f{r}{N} - \vep \ri ) \ri \rf  =
k$} and the lemma is true for $m = r$.  Similarly, {\small $\li \lf
n \li ( \f{r^\prime}{N} - \vep \ri ) \ri \rf  = k + 1$} and the
lemma is true for $m = r^\prime$.

Since $m$ is an integer, we have {\small $N \li (  \f{k}{n} + \vep
\ri ) \leq r < m \leq r^\prime - 1 < N \li ( \f{k+1}{n} + \vep \ri
)$} for $r < m < r^\prime$. Hence, $\f{k}{n} < \f{m}{N} - \vep <
\f{k+1}{n}$, or equivalently, $k < n \li (  \f{m}{N} - \vep \ri ) <
k + 1$, which implies that $\li \lf n \li (  \f{m}{N} - \vep \ri )
\ri \rf = k$ for $r < m < r^\prime$.

\epf

\beL \la{compa} Let $0 \leq \ro <  N$ and $g \leq h$. Then $S(n, g,
h + 1, \ro + 1, N) - S(n, g, h, \ro, N) \geq 0$. \eeL

\bpf

Note that, by Lemma \ref{dou}, \bee &   & S(n, g, h + 1, \ro + 1,
N)
- S(n, g, h, \ro, N)\\
 & = & \bi{\ro +1}{h+1} \bi{N - \ro - 1}{n - h - 1} \li \sh \bi{N}{n} \ri. +
 S(n, g, h, \ro + 1, N) - S(n, g, h, \ro, N)\\
& = &  \bi{\ro +1}{h+1} \bi{N - \ro - 1}{n - h - 1}
\li \sh \bi{N}{n} \ri. + T(g - 1, \ro, N, n) - T(h, \ro, N, n)\\
& = & \li [ \bi{\ro +1}{h+1} \bi{N - \ro - 1}{n - h - 1}  -
\bi{\ro}{h} \bi{N - \ro - 1}{n - h - 1} \ri ] \li \sh \bi{N}{n} \ri. + T(g - 1, \ro, N, n)\\
& = &  \bi{\ro}{h+1} \bi{N - \ro - 1}{n - h - 1} \li \sh \bi{N}{n}
\ri. + T(g - 1, \ro, N, n) \geq  0, \eee where the last equality
follows from (\ref{basic}).

 \epf

\beL \la{compb} Let $0 < \tau \leq N$ and $g \leq h$. Then, $S(n, g
- 1, h, \tau - 1, N) - S(n, g, h, \tau, N) \geq 0$. \eeL

\bpf

Note that, by Lemma \ref{dou}, \bee &   & S(n, g - 1, h, \tau - 1, N) - S(n, g, h, \tau, N) \\
& = & \bi{ \tau - 1 } {g - 1  } \bi{ N - \tau + 1  } { n - g + 1 }
\li \sh \bi{N}{n} \ri.
+ S(n, g, h, \tau - 1, N) - S(n, g, h, \tau, N)\\
& = &  \bi{ \tau - 1 } {g - 1  } \bi{ N - \tau + 1  } { n - g + 1
} \li \sh \bi{N}{n} \ri. + T(h, \tau - 1, N, n) - T(g - 1, \tau - 1, N, n)\\
& = &  \li [ \bi{ \tau - 1 } {g - 1  } \bi{ N - \tau + 1  } { n -
g + 1 } - \bi{ \tau - 1 } {g - 1  } \bi{ N
- \tau } { n - g } \ri ] \li \sh \bi{N}{n} \ri. + T(h, \tau - 1, N, n)\\
& = & \bi{ \tau - 1 } {g - 1  } \bi{ N - \tau } { n - g + 1 } \li
\sh \bi{N}{n} \ri. + T(h, \tau - 1, N, n) \geq 0,
 \eee
 where the last equality
follows from (\ref{basic}).

\epf

\beL

\la{last}

Let $\ro < \tau$ be two consecutive elements of the ascending
arrangement of all distinct elements of

{\small $\{L, U \} \cup \li \{ \li \lf N \li ( \f{k}{n} - \vep \ri
) \ri \rf \in (L, U): k \in \bb{Z} \ri \} \bigcup \li \{ \li \lc N
\li ( \f{k}{n} + \vep \ri ) \ri \rc \in (L, U): k \in \bb{Z} \ri
\}$. }

Then,
\[
C(M) = S(n, g(M), h(M), M,N) \geq \min \{ S(n, g(\ro), h(\ro), \ro,N), S(n, g(\tau), h(\tau), \tau,N) \}
\]
for $\ro \leq M \leq \tau$.

 \eeL

\bpf

Since the lemma is obviously true if $\tau = \ro + 1$, we may focus
on the situation that $\ro < \tau - 1$. To show the lemma, it
suffices to consider $7$ cases as follows.

\bsk

Case (i): $\ro = \li \lf N \li (  \f{k}{n} - \vep \ri ) \ri \rf$
and $\tau = \li \lc N \li (  \f{r}{n} + \vep \ri ) \ri \rc$. Note
that {\small \[ \li \lc N \li (  \f{r - 1}{n} + \vep \ri ) \ri \rc
\leq \li \lf N \li ( \f{k}{n} - \vep \ri ) \ri \rf = \ro < \tau =
\li \lc N \li ( \f{r}{n} + \vep \ri ) \ri \rc \leq \li \lf N \li (
\f{k+1}{n} - \vep \ri ) \ri \rf.
\]}
By Lemma \ref{lemg}, $g(\tau -1) = g(\tau) - 1 = r$. Since $\ro <
\tau - 1$, by Lemma \ref{lemh}, we have $h(\tau - 1) = h(\tau ) =
k$. By Lemma \ref{compb}, we have $S(n, g - 1, h, \tau - 1, N) -
S(n, g, h, \tau, N) \geq 0$.  That is, $C(\tau - 1) \geq C(\tau)$.

By Lemma \ref{lemh}, $h(\ro + 1) = h(\ro) + 1 = k$.  Since $\ro <
\tau - 1$, by Lemma \ref{lemg}, we have $g(\ro) = g(\ro + 1) = r$.
By Lemma \ref{compa}, we have $S(n, g, h + 1, \ro + 1, N) - S(n, g,
h, \ro, N) \geq 0$.  That is, $C(\ro + 1) \geq C(\ro)$.

By Lemma \ref{lem77}, we have $C(M) \geq \min \{ C(\ro + 1),\;
C(\tau -1) \}$ for $\ro + 1 \leq M \leq \tau - 1$.  It follows that
$C(M) \geq \min \{ C(\ro),\; C(\tau) \}$ for $\ro \leq M \leq \tau$.

\bsk

Case (ii): $\ro = \li \lf N \li (  \f{k}{n} - \vep \ri ) \ri \rf$
and $\tau = \li \lf N \li (  \f{k+1}{n} - \vep \ri ) \ri \rf$. In
this case, there exists an integer $l$ such that {\small \[ \li \lc
N \li ( \f{l}{n} + \vep \ri ) \ri \rc \leq \li \lf N \li ( \f{k}{n}
- \vep \ri ) \ri \rf = \ro < \tau = \li \lf N \li ( \f{k+1}{n} -
\vep \ri ) \ri \rf \leq \li \lc N \li ( \f{l+1}{n} + \vep \ri ) \ri
\rc.
\]}
By Lemma \ref{lemh}, we have $h(\ro) = k -1$ and $h(\ro + 1) = k$.
Since $\ro < \tau - 1$, by Lemma \ref{lemg}, we have $g(\ro+1) =
g(\ro) = l + 1$.  By Lemma \ref{compa}, $S(n, g, h + 1, \ro + 1, N)
- S(n, g, h, \ro, N) \geq 0$.  That is, $C(\ro + 1) \geq C(\ro)$.

For the value of $\tau$, there are two sub-cases:

(a) $\tau < \li \lc N \li ( \f{l+1}{n} + \vep \ri ) \ri \rc$. In
this case, by Lemma \ref{lemh} and Lemma \ref{lemg}, we have $g(M) =
g(\tau)$ and $h(M) = h(\tau)$ for $\ro < M \leq \tau$. It follows
from Lemma \ref{lem77} that $C(M) = S(n, g(M), h(M), M,N) \geq \min
\{ C(\ro +1), \; C(\tau) \}$ for $\ro < M \leq \tau$.  As a result,
$C(M) \geq \min \{ C(\ro),\; C(\tau) \}$ for $\ro \leq M \leq \tau$.

(b) $\tau = \li \lc N \li ( \f{l+1}{n} + \vep \ri ) \ri \rc$. This
case is the same as Case (i) previously studied.

\bsk

Case (iii): $\ro = \li \lc N \li (  \f{k}{n} + \vep \ri ) \ri \rc$
and $\tau = \li \lc N \li (  \f{k+1}{n} + \vep \ri ) \ri \rc$. In
this case, there exists an integer $l$ such that {\small \[ \li \lf
N \li ( \f{l}{n} - \vep \ri ) \ri \rf \leq \li \lc N \li ( \f{k}{n}
+ \vep \ri ) \ri \rc = \ro < \tau = \li \lc N \li ( \f{k+1}{n} +
\vep \ri ) \ri \rc \leq \li \lf N \li ( \f{l+1}{n} -\vep \ri ) \ri
\rf.
\]}
By Lemma \ref{lemg}, $g(\tau) = k + 2$ and $g(\tau - 1) = k + 1$.
Since $\ro < \tau - 1$, by Lemma \ref{lemh}, we have $h(\tau - 1) =
h(\tau) = l$. Thus, by Lemma \ref{compb}, we have $S(n, g - 1, h,
\tau - 1, N) - S(n, g, h, \tau, N) \geq 0$.  That is, $C(\tau - 1)
\geq C(\tau)$.

For the value of $\ro$, there are two sub-cases.

(a): $\ro > \li \lf N \li ( \f{l}{n} - \vep \ri ) \ri \rf$.
 Since $\ro < \tau - 1$, by Lemma \ref{lemg} and Lemma \ref{lemh}, we have $g(\ro+1) = g(\ro) = k
+1$ and $h(\ro+1) = h(\ro) = l$.  Hence, by Lemma \ref{lem77}, $C(M)
= S(n, g(M), h(M), M,N) \geq \min \{ C(\ro), \; C(\tau - 1) \}$ for
$\ro$ to $\tau - 1$.

(b): $\ro = \li \lf N \li ( \f{l}{n} - \vep \ri ) \ri \rf$. This
case becomes Case (i) previously studied.

\bsk

 Case (iv): $\ro = \li \lc N \li (  \f{k}{n} +
\vep \ri ) \ri \rc$ and $\tau = \li \lf N \li (  \f{r}{n} - \vep
\ri ) \ri \rf$. Note that {\small \[ \li \lf N \li (  \f{r - 1}{n}
- \vep \ri ) \ri \rf \leq \li \lc N \li ( \f{k}{n} + \vep \ri )
\ri \rc = \ro < \tau = \li \lf N \li ( \f{r}{n} - \vep \ri ) \ri
\rf \leq \li \lc N \li ( \f{k+1}{n} + \vep \ri ) \ri \rc.
\]}
If {\small $\li \lf N \li (  \f{r - 1}{n} - \vep \ri ) \ri \rf = \li
\lc N \li ( \f{k}{n} + \vep \ri ) \ri \rc$}, then the case becomes
{\small $\li \lf N \li (  \f{r - 1}{n} - \vep \ri ) \ri \rf = \ro <
\tau = \li \lf N \li ( \f{r}{n} - \vep \ri ) \ri \rf$},  which has
been studied previously in Case (ii).

If {\small $\li \lf N \li ( \f{r}{n} - \vep \ri ) \ri \rf = \li \lc
N \li ( \f{k+1}{n} + \vep \ri ) \ri \rc$}, then the case becomes
{\small $\li \lc N \li ( \f{k}{n} + \vep \ri ) \ri \rc = \ro < \tau
= \li \lc N \li ( \f{k+1}{n} + \vep \ri ) \ri \rc$, } which has been
studied previously in Case (iii). So, we only need to consider the
situation that {\small \[ \li \lf N \li (  \f{r - 1}{n} - \vep \ri )
\ri \rf < \li \lc N \li ( \f{k}{n} + \vep \ri ) \ri \rc = \ro < \tau
= \li \lf N \li ( \f{r}{n} - \vep \ri ) \ri \rf < \li \lc N \li (
\f{k+1}{n} + \vep \ri ) \ri \rc.
\]}
Since $\ro < \tau - 1$, by Lemma \ref{lemh} and Lemma \ref{lemg}, we
have
\[
g(\ro) =  g(\ro + 1) = k + 1, \qqu h(\ro) =  h(\ro + 1) = r - 1,
\]
\[
h(\tau - 1) = h(\tau) = r - 1, \qqu g(\tau - 1) = g(\tau) = k + 1.
\]
It follows that $g(M) = k + 1$ and $h(M) = r - 1$ for $\ro \leq M
\leq \tau$. By Lemma \ref{lem77}, $C(M) \geq \min \{ C(\ro), \;
C(\tau) \}$ for $\ro \leq M \leq \tau$.

\bsk

Case (v):  $\ro = L$ and $\tau \in \mscr{Q}$  with \[ \mscr{Q} = \li
\{ \li \lf  N \li ( \f{k}{n} - \vep \ri ) \ri \rf : k \in \bb{Z} \ri
\}
 \bigcup \li \{ \li \lc N \li ( \f{k}{n} + \vep \ri )
\ri \rc : k \in \bb{Z} \ri \}.
\]
The case that $\ro = L \in \mscr{Q}$ can be included in previous cases. So we can focus on the case that $\ro
\notin \mscr{Q}$.  In this case, by Lemma \ref{lem77}, $C(M) \geq \min \{ C(\ro), \; C(\tau - 1) \}$ for $\ro
\leq M \leq \tau - 1$. The comparison of $C(\tau - 1)$ with $C(\tau)$ is as before.

\bsk

Case (vi):  $\ro \in \mscr{Q}$ and $\tau = U$. The case that $\tau = U \in \mscr{Q}$ can be included in previous
cases. So we can focus on the case that $\tau \notin \mscr{Q}$.
 In this case, by Lemma \ref{lem77}, $C(M) \geq \min \{ C(\ro + 1), \;
C(\tau) \}$ for $\ro + 1 \leq M \leq \tau$. The comparison of $C(\ro
+ 1)$ with $C(\ro)$ is as before.

\bsk

Case (vii): $\ro = L$ and $\tau = U$.  The cases that $L$ or $U$ belongs to $\mscr{Q}$ can be considered as Case
(v) or Case (vi). Thus, we can focus on the case that neither $L$ nor $U$ belongs to $\mscr{Q}$.  Invoking Lemma
\ref{lem77},  we have $C(M) \geq \min \{ C(\ro), \; C(\tau) \}$ for $\ro \leq M \leq \tau$.

 \epf

 \bsk

 Now we are in position to prove Theorem \ref{thm_abs}.
 Clearly, the statement about the coverage probability follows immediately follows from Lemma \ref{last}.
It remains to compute the number of elements in the set.  Making use of the fact that, for any real number $x$
and integer $r$,
\[
\lf x \rf < r \LRA x < r, \qqu  \lf x \rf > r \LRA x \geq 1+  r,
 \]
we have
\[
\li \lf N \li ( \f{k}{n} - \vep \ri ) \ri \rf  \in (L, U) \LRA n
\li ( \f{L +1}{N} + \vep \ri ) \leq k < n \li ( \f{U}{N} + \vep
\ri ).
\]
That is,
\[
\li \lf N \li ( \f{k}{n} - \vep \ri ) \ri \rf  \in (L, U) \LRA \li
\lc n \li ( \f{L +1}{N} + \vep \ri ) \ri \rc \leq k \leq  \li \lc
n \li ( \f{U}{N} + \vep \ri ) \ri \rc - 1.
\]
Thus, the number of elements in the set $\li \{ \li \lf N \li ( \f{k}{n} - \vep \ri ) \ri \rf  \in (L, U): k \in
\bb{Z} \ri \}$ is {\small \[ \li \lc n \li ( \f{U}{N} + \vep \ri ) \ri \rc - \li \lc n \li ( \f{L+1}{N} + \vep
\ri ) \ri \rc < n \li ( \f{U}{N} + \vep \ri ) -  n \li ( \f{L+1}{N} + \vep \ri ) + 1 = 1 + n \li ( \f{ U - L - 1
} { N } \ri ),
\]}
where we have used inequalities $\lc x \rc - \lc y \rc < x - \lc y \rc + 1 \leq x - y + 1$ for real numbers $x$
and $y$.

\bsk

Similarly,  making use of the fact that, for any real number $x$ and integer $r$, \[ \lc x \rc < r \LRA x \leq r
- 1, \qqu \lc x \rc > r \LRA x > r,
 \]
we have
\[
\li \lc N \li ( \f{k}{n} + \vep \ri ) \ri \rc  \in (L, U) \LRA n
\li ( \f{L}{N} - \vep \ri ) < k \leq n \li ( \f{U - 1}{N} - \vep
\ri ).
\]
That is,
\[
\li \lc N \li ( \f{k}{n} + \vep \ri ) \ri \rc  \in (L, U) \LRA \li
\lf n \li ( \f{L}{N} - \vep \ri ) \ri \rf + 1 \leq k \leq  \li \lf
n \li ( \f{U - 1}{N} - \vep \ri ) \ri \rf.
\]
Thus, the number of elements in the set $\li \{ \li \lc N \li ( \f{k}{n} + \vep \ri ) \ri \rc  \in (L, U): k \in
\bb{Z} \ri \}$ is {\small \[ \li \lf n \li ( \f{U - 1}{N} - \vep \ri ) \ri \rf - \li \lf n \li ( \f{L}{N} - \vep
\ri ) \ri \rf < n \li ( \f{U - 1}{N} - \vep \ri ) -  n \li ( \f{L}{N} - \vep \ri ) + 1 = 1 + n \li ( \f{ U - L -
1 } { N } \ri ),
\]}
where we have used inequalities $\lf x \rf - \lf y \rf \leq x -
\lf y \rf < x - y + 1$.

\bsk

Therefore, the total number of elements in

$\li \{ \li \lf N \li ( \f{k}{n} - \vep \ri ) \ri \rf  \in (L, U): k \in \bb{Z} \ri \} \cup \li \{ \li \lc N \li
( \f{k}{n} + \vep \ri ) \ri \rc  \in (L, U): k \in \bb{Z} \ri \}$

 is less than
\[
1 + n \li ( \f{ U - L - 1 } { N } \ri ) + 1 + n \li ( \f{ U - L -
1 } { N } \ri ) = 2 + 2 n \li ( \f{ U - L - 1} { N } \ri ).
\]
The bound becomes $4 + 2 n \li ( \f{ U - L - 1} { N } \ri )$ since $L$ and $U$ are needed to be counted.  This
concludes the proof of Theorem \ref{thm_abs}.

\end{document}